\documentclass[final,onefignum,onetabnum,letterpaper]{siamonline190516}

\usepackage{xr-hyper} 
\usepackage{xr}

\usepackage{hyperref}

\usepackage[utf8]{inputenc}
\usepackage{url}
\usepackage{soul}
\usepackage{lipsum}
\usepackage{amsfonts}
\usepackage{epstopdf}
\ifpdf
  \DeclareGraphicsExtensions{.eps,.pdf,.png,.jpg}
\else
  \DeclareGraphicsExtensions{.eps}
\fi

\usepackage{datetime}
\newdateformat{monthyeardate}{%
  \monthname[\THEMONTH] \THEDAY, \THEYEAR}

\usepackage{academicons}
\usepackage{xcolor}

\usepackage{amsmath}
\allowdisplaybreaks
\usepackage{amssymb}
\usepackage{hyperref,cleveref}
\usepackage{commath}
\usepackage{mathtools}
\usepackage{bbm}

\usepackage{color}
\usepackage{graphicx}
\usepackage[small]{caption}
\usepackage{subcaption}

\usepackage{relsize}
\usepackage{adjustbox}
\usepackage{algorithm,algpseudocode}
\usepackage{booktabs}
\usepackage{tikz}

\usepackage{verbatim}

\usepackage{enumitem}
\setlist[enumerate]{leftmargin=.5in}
\setlist[itemize]{leftmargin=.5in}

\newsiamthm{problem}{Problem}
\newsiamremark{remark}{Remark}
\newsiamremark{hypothesis}{Hypothesis} 
\crefname{hypothesis}{Hypothesis}{Hypotheses}
\newsiamremark{example}{Example}
\newsiamthm{claim}{Claim}
\newsiamthm{conjecture}{Conjecture}

\headers{Hybrid ABBA GMRES}{R.\ Bentley, M.\ Pasha, M. Sabaté Landman, L.\ Yang, J.\ Zhang}

\title{Hybrid ABBA-GMRES for Unmatched Backprojectors in Large Scale X-Ray Computerized Tomography
\thanks{
\monthyeardate\today 
\corresponding{Mirjeta Pasha} 
}}

\author{
Ryan Bentley\footnotemark[2]
\and 
Mirjeta Pasha\thanks{Department of Mathematics \& Academy of Data Science, Virginia Tech, USA (\email{rbentley5@vt.edu}, \email{mpasha@vt.edu}, \email{jefferyz@vt.edu})}
\and
Malena Sabaté Landman\thanks{
Department of Mathematical Sciences, University of Bath, UK (\email{msl39@bath.ac.uk})}
\and 
Luisa Yang\thanks{
Mathematical Institute, University of Oxford, UK (\email{luisayang27@gmail.com})}
\and
Jeffery Zhang \footnotemark[2]
}

\renewcommand{\t} {^{\top}}                                
\usepackage{amsopn}

\usepackage{ifthen}
\usepackage{float}
\floatstyle{ruled}
\newfloat{algorithm}{htb}{alg}
\floatname{algorithm}{Algorithm}
\newcounter{algo@row}
\newcounter{algo@rowindent}
\newcommand{\algofont}[1]{\textbf{#1}}
\newcommand{\algonumbersize}[1]{\scriptsize{#1}}
\newcommand{\algopreitem}[1][\arabic{algo@row}]{\texttt{\algonumbersize{#1}}}
\newcommand{\algoitemskip}{\hspace{\value{algo@rowindent}cc}}
\newcommand{\algonewnestedopen}[2]{
	\newcommand{#1}[1][]{%
		\ifthenelse{\equal{##1}{}}{\item}{\item[{\algopreitem[##1]}]}
		\algoitemskip\algofont{#2}%
		\addtocounter{algo@rowindent}{1}%
		\ignorespaces
	}
}
\newcommand{\algonewnestedaux}[2]{
	\newcommand{#1}[1][]{
		\addtocounter{algo@rowindent}{-1}
		\ifthenelse{\equal{##1}{}}{\item}{\item[{\algopreitem[##1]}]}
		\algoitemskip\algofont{#2}%
		\addtocounter{algo@rowindent}{+1}%
		\ignorespaces
	}
}
\newcommand{\algonewnestedclose}[2]{
	\newcommand{#1}[1][]{
		\addtocounter{algo@rowindent}{-1}
		\ifthenelse{\equal{##1}{}}{\item}{\item[{\algopreitem[##1]}]}
		\algoitemskip\algofont{#2}%
		\ignorespaces
	}
}
\newcommand{\algonewcommand}[2]{
	\newcommand{#1}[1][default]{
		\ifthenelse{\equal{##1}{default}}{\item}{\item[{\algopreitem[##1]}]}%
		\algoitemskip\algofont{#2}%
		\ignorespaces
	}%
}
\newcommand{\algonewkeyword}[2]{\newcommand{#1}{\algofont{#2}}}
\algonewcommand{\STATE}{\ignorespaces}
\algonewcommand{\INPUT}{Input: }
\algonewcommand{\pINPUT}{\phantom{Input: }}
\algonewcommand{\COMPUTE}{Compute: }
\algonewcommand{\OUTPUT}{Output: }
\algonewcommand{\pOUTPUT}{\phantom{Output: }}

\algonewnestedopen{\IF}{if }
\algonewnestedaux{\ELSEIF}{else if }
\algonewnestedaux{\ELSE}{else }
\algonewnestedclose{\ENDIF}{end if }
\algonewnestedopen{\FOR}{for }
\algonewnestedclose{\ENDFOR}{end for }
\algonewnestedopen{\WHILE}{while }
\algonewnestedclose{\ENDWHILE}{end while }
\algonewcommand{\BREAK}{break}%

\algonewkeyword{\To}{to }%
\algonewkeyword{\Do}{do }%
\algonewkeyword{\Then}{then }%
\algonewkeyword{\End}{end }%
\algonewkeyword{\AND}{and }%
\algonewkeyword{\True}{true }%
\algonewkeyword{\False}{false }%
\algonewkeyword{\irbleigs}{irbleigs }%
\algonewkeyword{\tridiag}{tridiag}%
\algonewkeyword{\reorth}{reorth}%

\DeclareMathOperator*{\argmin}{arg\,min}

\newcommand{\R}{\mathbb{R}}

\newcommand{\bA}{{\bf A}}
\newcommand{\bB}{{\bf B}}

\newcommand{\bH}{{\bf H}}
\newcommand{\bI}{{\bf I}}

\newcommand{\bM}{{\bf M}}

\newcommand{\bW}{{\bf W}}

\newcommand{\bb}{{\bf b}}

\newcommand{\be}{{\bf e}}

\newcommand{\bh}{{\bf h}}

\newcommand{\bq}{{\bf q}}
\newcommand{\br}{{\bf r}}

\newcommand{\bu}{{\bf u}}

\newcommand{\bw}{{\bf w}}
\newcommand{\bx}{{\bf x}}
\newcommand{\by}{{\bf y}}
\newcommand{\bz}{{\bf z}}
\newcommand{\bzero}{{\bf0}}

\makeatletter
\newcommand*{\addFileDependency}[1]{
  \typeout{(#1)}
  \@addtofilelist{#1}
  \IfFileExists{#1}{}{\typeout{No file #1.}}
}
\makeatother

\usepackage{lineno}
\linenumbers

\ifpdf
\hypersetup{
  pdftitle={Priorconditioned projection methods},
  pdfauthor={Pasha}
}
\fi
\begin{document}
\nolinenumbers
\maketitle
\begin{abstract} In large-scale X-ray computed tomography (CT), matrix-free iterative methods are essential due to the prohibitive cost of explicitly forming the system matrix. In practice, forward projectors and backprojectors are often implemented with different discretizations or accelerations, leading to unmatched projector pairs. This mismatch violates the adjointness assumptions underlying classical least-squares solvers, so the resulting iterations no longer correspond to a true least-squares problem and can exhibit non-symmetric or inconsistent behavior.
Prior work has explored Krylov subspace solvers such as AB-GMRES and BA-GMRES to handle unmatched projector pairs, where these methods exhibit semi-convergent regularizing behavior. Under matched conditions, AB-GMRES and BA-GMRES reduce to LSQR and LSMR, respectively. However, in the presence of unmatched projectors, AB- and BA-GMRES have been observed to yield improved reconstruction quality compared to classical least-squares solvers. In this paper, we develop hybrid AB- and BA-GMRES methods that incorporate Tikhonov regularization directly into the Krylov subspace iterations.
We also examine the relationship between the proposed methods and hybrid variants of LSQR and LSMR, considering both matched and unmatched backprojectors. We propose automatic strategies for selecting regularization parameters, including approaches based on the L-curve and generalized cross validation (GCV), and analyze their effect on convergence behavior and image quality. Numerical experiments on two-dimensional CT problems using GPU-accelerated projectors demonstrate that the proposed hybrid AB- and BA-GMRES methods mitigate semi-convergence, produce higher-quality reconstructions, and exhibit more stable stopping behavior than their non-hybrid counterparts.
\end{abstract}


\section{Introduction}\label{sec: Intro}
X-ray computed tomography reconstructs an image of an object’s internal attenuation properties from line integrals measured along many rays. After discretizing the imaging geometry, the reconstruction problem takes the form
\begin{equation}\label{eq:linear_system}
\bA\bx \approx \bb,
\end{equation}
where $\bA\in\mathbb{R}^{m\times n}$ encodes the forward projection, $\bx\in\mathbb{R}^n$ is the discretized image, and $\bb~\in~\mathbb{R}^m$ represents the measured data, which includes additive measurement noise. In large-scale imaging applications, $\bA$ is too large to store explicitly; instead, matrix-vector products with $\bA$ and its adjoint are computed using optimized projector and backprojector routines. In practice, these routines often arise from differing discretizations or hardware-specific approximations, leading to \emph{unmatched projector pairs} in which the backprojector $\bB$ is not the exact adjoint $\bA^\mathsf{T}$ \cite{zeng2000effects,hansen2022gmres}. For example, unmatched projector/backprojector pairs commonly appear in modern iterative CT implementations because approximate backprojectors can significantly reduce computation time and enhance hardware utilization. 

However, mismatched operators result in non-symmetric normal equations, and classical iterative least-squares methods such as CGLS or other gradient-based methods typically fail to deliver meaningful reconstructions of $\bx$ in \eqref{eq:linear_system} \cite{zeng2000effects,sidky2022gmres}. Early work on this issue pointed out fundamental constraints on unmatched pairs and analyzed the spectral properties that influence convergence \cite{zeng2000effects}. To address the mismatch directly, preconditioned versions of the GMRES algorithm — specifically the AB- and BA- GMRES methods —  were proposed \cite{hansen2022gmres}. These solvers form square systems from unmatched forward/backprojector pairs (namely $\bA\bB$ and $\bB\bA$) and apply GMRES to the resulting unmatched normal equations without requiring an explicit adjoint, and can serve as regularizing iterative solvers. 
However, one of the drawbacks of AB-/BA-GMRES methods is that they suffer from \emph{semi-convergence}: the reconstruction error typically decreases initially but deteriorates after a few iterations due to the ill-posedness of the underlying inverse problem \cite{hansen2022gmres, hansen1998rank}. 

Semi-convergence is a common challenge of iterative methods for linear discrete inverse problems, and has motivated the development of \emph{hybrid regularization} strategies, that embed Tikhonov regularization into Krylov subspace methods \cite{hansen1998rank,chung2024survey}. Such hybrid methods automatically balance data fidelity and regularization within the iterative process and have proven effective for large-scale inverse problems. 

In this work, we present novel hybrid AB- and BA-GMRES methods specifically designed for settings with unmatched backprojectors. By incorporating explicit regularization into the projected problem, the proposed hybrid methods mitigate the semi-convergence observed in their standard counterparts and yield reconstructions that are more robust to noise, operator mismatch, and modeling errors.
Moreover, because the regularization is applied in a low-dimensional Krylov subspace, the additional computational cost is negligible compared to the cost of forward and backprojection operations. The numerical experiments illustrate that hybrid AB-/BA-GMRES methods consistently mitigate semi-convergence, produce improved reconstructions, and offer stable stopping behavior. 
Our findings support the adoption of hybrid AB-/BA-GMRES frameworks for regularized CT reconstruction with unmatched backprojectors, and provide practical guidance for implementation in large-scale imaging systems.

Our contributions can be summarized as follows:
\begin{itemize}
\item We formulate hybrid AB-/BA-GMRES variants that incorporate Tikhonov regularization within the Krylov iteration for unmatched projector pairs. We analyze the relationship between these methods and  hybrid LSQR/LSMR in the matched projector case, i.e. when $\bB=\bA^\mathsf{T}$; as well as showing that ‘project-then-regularize’ is not equivalent to ‘regularize-then-project’ \cite{hansen1998rank} for AB-/BA-GMRES methods.

\item We investigate automatic strategies for selecting regularization parameters at each iteration, including the L-curve and generalized cross validation (GCV).

\item We provide numerical experiments comparing AB-/BA-GMRES and LSQR/LSMR with and without explicit regularization in the projected problem, both in matched and unmatched settings. Moreover, we provide a thorough numerical evaluation across different CT problems using GPU-accelerated projectors.
\end{itemize}

The paper is organized as follows. Section \ref{sec:background} presents the mathematical framework for solving linear discrete inverse problems with unmatched backprojectors, existing Krylov subspace methods for this setting, and hybrids regularization methods. Section \ref{sec:hybridABBA} introduces the new method, giving theoretical considerations and analyzing regularization parameter choices in Section \ref{sec:regparam} as well as stopping criteria in  Section \ref{sec:stop_critieria}. Finally, numerical experiments are presented in Section \ref{sec:computational}. 

\section{Background and preliminaries}\label{sec:background} 
The reconstruction problem \eqref{eq:linear_system} is commonly posed as a least-squares problem. When $m \geq n$, this is equivalent to solving the normal equations
\begin{equation}
\bA^{\mathsf{T}} \bA \bx = \bA^{\mathsf{T}} \bb,
\label{eq:matched_normal_1}
\end{equation}
while the minimum norm solution of the least-squares problem for $m\leq n$ satisfies
\begin{equation}
 \bA \bA^{\mathsf{T}} \bu =  \bb, \quad \bx = \bA^{\mathsf{T}}\bu.
\label{eq:matched_normal_2}
\end{equation}

In practical implementations, the transpose $\bA^{\mathsf{T}}$ is realized numerically through a backprojection operator which we denote by $\bB \in \mathbb{R}^{n \times m}$, so that \eqref{eq:matched_normal_1} is implicitly replaced by 
\begin{equation}
\bB \bA \bx = \bB \bb,
\label{eq:normal_general_1}
\end{equation}
and \eqref{eq:matched_normal_2} is effectively replaced by
\begin{equation}
\bA \bB \bu =  \bb , \quad \bx = \bB \bu.
\label{eq:normal_general_2}
\end{equation}
When $\bB = \bA^{\mathsf{T}}$, \eqref{eq:normal_general_1} corresponds to the \emph{matched} normal equations, and both  \eqref{eq:normal_general_1} and  \eqref{eq:normal_general_2} can be solved using classical iterative solvers such as CGLS or the mathematically equivalent LSQR. 
In this paper, we focus on to \emph{unmatched projector pairs} for which
\begin{equation*}
\bB \neq \bA^{\mathsf{T}}.
\end{equation*}
Such unmatched projectors are routinely encountered in modern CT software frameworks, where optimized GPU implementations frequently employ approximate backprojectors to reduce memory usage and computational cost, e.g. the ASTRA Toolbox \cite{ASTRA}, or the TIGRE toolbox \cite{biguri2016tigre,biguri2025tigre}.  
Note that, when the range of $\bB$ is the same than the range of $\bA\t$, equation~\eqref{eq:normal_general_1} is called the not-normal equations~\cite{wathen2025nonnormal}, and have the same solution as \eqref{eq:matched_normal_1}.

To address the unmatched setting directly, the AB-/BA-GMRES methods apply GMRES to square systems formed from unmatched forward and backprojectors \cite{hansen2022gmres}. These methods avoid the need for an explicit adjoint and provide a principled framework for reconstruction with unmatched operators. However, as with other iterative methods for ill-posed inverse problems, AB-/BA-GMRES exhibit semi-convergence: early iterations act as a form of implicit regularization, while later iterations amplify noise.

The development of robust and efficient solvers for unmatched projector pairs therefore remains an important challenge in large-scale CT reconstruction. In the following sections, we build upon the AB-/BA-GMRES framework and introduce hybrid regularization strategies designed to mitigate semi-convergence and improve reconstruction stability in the presence of unmatched backprojectors.
\subsection{ABBA Iterative Methods}
\label{sec:ABBA}

The AB-GMRES and BA-GMRES methods were introduced in \cite{hansen2022gmres} and further analyzed in \cite{Knudsen_2023} as Krylov subspace solvers for large-scale tomographic reconstruction problems with unmatched projector pairs. These methods are specifically designed to handle rectangular systems arising in computed tomography in this case, when the matrices appearing in equations \eqref{eq:normal_general_1} and \eqref{eq:normal_general_2} are no longer symmetric, and conjugate gradient–based methods are not applicable. The ABBA framework circumvents this difficulty by applying GMRES to square systems formed from compositions of $\bA$ and $\bB$, without requiring an explicit adjoint.

\paragraph{AB-GMRES}
AB-GMRES is primarily intended for underdetermined systems ($m < n$). The method applies GMRES to the square system \eqref{eq:normal_general_2}, where operator $\bA\bB \in \mathbb{R}^{m \times m}$ is square but generally non-symmetric, and can be easily implemented in a matrix-free fashion, only requiring evaluations of the action of $\bA$ and $\bB$. In particular, at each iteration $k$, AB-GMRES finds the solution
\begin{equation}
    \bu_k =\arg\min_{\bu\in\mathcal{K}_k(\bA\bB,\bb)}\|\bA\bB\bu-\bb\|_2^2, \quad \bx_k = \bB \bu_k \in \bx_0+ \mathcal{K}_k(\bB\bA,\bB\bb).
    \label{eq:AB_optimality}
\end{equation}

\paragraph{BA-GMRES}
BA-GMRES is more naturally suited for overdetermined systems ($m \geq n$). In this case, GMRES is applied to the system \eqref{eq:normal_general_1}.  In particular, at each iteration $k$, BA-GMRES find the solution
\begin{equation}
    \bx_k =\arg\min_{\bx\in\mathcal{K}_k(\bB\bA,\bB\bb)}\|\bB\bA\bx-\bB\bb\|_2^2 \in \bx_0+ \mathcal{K}_k(\bB\bA,\bB\bb).
    \label{eq:BA_optimality}
\end{equation}
\paragraph{Algorithmic structure}
Both AB-GMRES and BA-GMRES are based on the Arnoldi process and construct Krylov subspaces associated with $\bA\bB$ or $\bB\bA$, respectively. At iteration $k$, the solution is obtained by solving a small projected least-squares problem involving the Hessenberg matrix produced by the Arnoldi iteration. Consider $\bM$ to be either $\bA\bB$ or $\bB\bA$, then the Arnoldi method produced an orthonormal Krylov basis $\bW_k = [\bw_1,\ldots,\bw_k]$, where $\bw_k$ is proportional to $\br_0$, and an associated upper Hessenberg matrix $\bH_k$ such that
\begin{equation}\label{eq:Hessenberg}
\bM\bW_k = \bW_{k+1} \bH_k.
\end{equation}
Note that, given an initial guess for the solution $\bx_0$, the initial residuals are
$\br_0 = \bb - \bA \bx_0$ for AB-GMRES and $\br_0 = \bB\bb - \bB\bA \bx_0$ for BA-GMRES.
The factorization~\eqref{eq:Hessenberg} can then be applied to efficiently solve the
minimization problem in~\eqref{eq:AB_optimality} or~\eqref{eq:BA_optimality} by reducing it
to a projected problem of the form
\begin{equation}
\by_k = \arg\min_{\by} \bigl\| \bH_k \by - \|\br_0\|_2 \be_1 \bigr\|_2.
\label{eq:proj_pbm}
\end{equation}
The resulting approximation at iteration $k$ is given by
\begin{equation}
\bx_k =
\begin{cases}
\bx_0 + \bB \bW_k \by_k, & \text{for AB-GMRES}, \\[2mm]
\bx_0 + \bW_k \by_k, & \text{for BA-GMRES}.
\end{cases}
\end{equation}
In both cases, the computational core of the method consists of solving a small
least-squares problem whose dimension grows with the iteration count, as in standard
GMRES. This problem can be handled efficiently using Givens rotations or a QR
factorization of $\bH_k$. The distinction between AB-GMRES and BA-GMRES lies solely in the
reconstruction of the approximate solution, with AB-GMRES applying the operator $\bB$
explicitly in the update, while BA-GMRES incorporates $\bB$ implicitly through the
definition of the residual. This difference may influence numerical stability and the
effectiveness of preconditioning.

\paragraph{ABBA-GMRES Semi-Convergence}

Like many Krylov subspace solvers applied to ill-posed inverse problems, AB-/BA-GMRES methods exhibits \emph{semi-convergence}. During the initial iterations, the iterates $\bx_k$ rapidly approach the true solution $\bx_{\mathrm{true}}$, as the associated Krylov subspaces capture components corresponding to the dominant spectral modes of the underlying operator. As the iteration count increases, however, components associated with small singular values—typically dominated by measurement noise and modeling errors—enter the solution, leading to noise amplification and a subsequent increase in reconstruction error.

This behavior can be understood through the implicit filtering properties of Krylov subspace methods. Early GMRES iterates effectively act as low-pass filters: the solution $\bx_k$ can be interpreted as a filtered expansion in terms of the singular vectors of $\bA$ (or of the composite operators $\bA\bB$ and $\bB\bA$), where the filters strongly damp components associated with small singular values. As the Krylov subspace dimension grows, these filters become less selective, eventually allowing noise-dominated components to contaminate the reconstruction. Similar phenomena are well documented for LSQR and CGLS in the matched setting \cite{hansen2010discrete}, and have been observed for GMRES applied to non-symmetric perturbations of the normal equations~\cite{hansen2022gmres,Knudsen_2023}.

In the ABBA framework, semi-convergence is further influenced by the use of unmatched projector/backprojector pairs. When $\bB \neq \bA^{\mathsf{T}}$, the composite operators $\bA\bB$ and $\bB\bA$ are generally non-normal, and their spectral properties may differ significantly from those of the matched normal equations. As a result, the ordering and separation of spectral components within the Krylov subspace can be less favorable, potentially accelerating the onset of semi-convergence or increasing sensitivity to noise. Despite this, AB- and BA-GMRES retain the characteristic regularizing behavior of Krylov methods during early iterations, making them viable iterative regularization schemes provided that appropriate safeguards are employed.

From a practical perspective, semi-convergence implies that the iteration count itself acts as a regularization parameter. However, standard stopping criteria based solely on residual norms may be inadequate, as the residual can continue to decrease even while the reconstruction error increases. Despite the fact that alternative stopping criteria exist, the reconstruction quality is usually very sensitive to this choice, so one might consider adding explicit regularization to obtain good quality reconstructions that are less sensitive to parameter tuning.

\subsection{Hybrid Methods}\label{sec:hybrid}
One of the most widely used approaches for stabilizing ill-posed inverse problems is Tikhonov regularization \cite{hansen2010discrete}. In this framework, the reconstruction is obtained by solving a penalized least-squares problem of the form
\begin{equation}
\min_{\bx \in \mathbb{R}^n}
\left\{
\|\bA\bx - \bb\|_2^2 + \lambda^2 \|\bx\|_2^2
\right\},
\label{eq:tikhonov}
\end{equation}
where $\lambda > 0$ is the regularization parameter. This balances the effect of the data fidelity term, which enforces consistency between the reconstructed image $\bx$ and the measured data $\bb$, and the Tikhonov regularization, which penalizes large solution norms and promotes smoothness by suppressing high-frequency components that are typically amplified by noise. 

From a spectral perspective, Tikhonov regularization acts as a filter on the singular values of $\bA$, damping contributions associated with small singular values that are most sensitive to noise. This explicit filtering complements the implicit regularization provided by early termination of Krylov subspace methods and provides a principled mechanism for controlling semi-convergence.

However, the regularization parameter $\lambda$ is not usually known a priori, and the selection of the regularization parameter $\lambda$ is critical to the success of Tikhonov regularization and remains one of its most challenging aspects. Common parameter-choice strategies include the discrepancy principle, the L-curve criterion, and generalized cross validation (GCV) \cite{hansen2010discrete}. 

In large-scale problems, finding a good parameter $\lambda$  might require solving \eqref{eq:tikhonov} multiple times, which can be very computationally demanding. A very attractive alternative in this context is the use of hybrid methods. These are a class of iterative projection regularization methods which, in its basic for, include Tikhonov regularization on the projected problem, see recent survey \cite{chung2024computational}. Since each of the projected problems is very small, $\lambda$ can be determined automatically and adaptively during the reconstruction process.

By incorporating the regularization term directly into the projected problems solved at each iteration, we propose new hybrid AB-/BA-GMRES methods, which provide improved stability, mitigate semi-convergence, and reduce reliance on ad hoc stopping criteria. These hybrid formulations are developed in the following section.
\section{Hybrid ABBA-GMRES}
\label{sec:hybridABBA} In this section, we propose new \emph{hybrid AB-GMRES}  and \emph{hybrid BA-GMRES} methods, which incorporate Tikhonov regularization explicitly in the projected Krylov subspace iterations. As mentioned in Section \ref{sec:hybrid}, the key idea of hybrid methods is to replace the projected least-squares problem solved at each AB-/BA-GMRES iteration, \eqref{eq:AB_optimality} or \eqref{eq:BA_optimality}, with a regularized problem defined on the current Krylov subspace. More explicitly, rather than solving \eqref{eq:proj_pbm}, 
we instead solve a Tikhonov-regularized problem of the form
\begin{equation}
\by^\text{h}_k=\arg\min_{\by}
\left\|
\begin{pmatrix}
\bH_k \\
\lambda \bI
\end{pmatrix}
\by -
\|\br_0\|_2
\begin{pmatrix}
\be_1 \\
\bzero
\end{pmatrix}
\right\|_2^2,
\label{eq:hybrid_proj}
\end{equation}
where $\bH_k$ is the Hessenberg matrix generated by the Arnoldi process applied to either $\bA\bB$ (AB-GMRES) or $\bB\bA$ (BA-GMRES), and $\lambda > 0$ is the regularization parameter.

Note that this formulation does not corresponds to applying Tikhonov regularization to the original inverse problem, but only to the restricted  $k$-dimensional Krylov subspace.
The resulting iterate
\begin{equation}
\bx^{\text{h,AB-/BA-GMRES}}_k =
\begin{cases}
\bx_0 + \bB \bW_k \by^\text{h}_k, & \text{hybrid AB-GMRES}, \\
\bx_0 + \bW_k \by^\text{h}_k, & \text{hybrid BA-GMRES},
\end{cases}
\end{equation}
balances data fidelity and solution smoothness within the Krylov subspace, thereby providing explicit control over noise amplification at each iteration. The resulting algorithms are summarized in Algorithms~\ref{alg:bagmres1} and~\ref{alg:bagmres2}.

\paragraph{Optimality conditions} Although the methods will not be run to convergence in practice, it is important to understand the limiting solution when the number of iterations reaches $\min(m,n)$. At iteration $n$, assuming $m \geq n$ and that not break-down of the algorithm has happened, the BA-GMRES method will produce the following solution
\begin{eqnarray*}
  \bx_n^{\text{h,BA-GMRES}}  &=& ((\bB\bA)^\mathsf{T}\bB\bA + \lambda \bI_n)^{-1}(\bB\bA)^\mathsf{T}\bB\bb \\
  &=& \argmin_{\bx}\|\bB\bA\bx-\bB\bb\|_2^2+\lambda^2\|\bx\|,
\end{eqnarray*}
where $\bI_n$ is the identity of dimension $n$. However, at iteration $m$, assuming $n \geq m$ and that not break-down of the algorithm has happened, the AB-GMRES method will produce the following solution
\begin{equation*}
  \bx_m^{\text{h,AB-GMRES}} = \bB\bz_m^{\text{h,AB-GMRES}} = \bB( ({\bA\bB})^\mathsf{T} {\bA\bB} + \lambda\bI_m)^{-1}{(\bA\bB)^\mathsf{T}}{\bb}.
\end{equation*}
Note that both hybrid BA-GMRES methods are `project-then-regularize' methods, and they are not equivalent to a `regularize-then-project' scheme. In this case, hybrid regularization needs to be understood as a way of stabilizing the iterations and avoiding semi-convergence, making the method less susceptible to the stopping criteria. However, these methods are not equivalent to projecting a Tikhonov problem, so the regularization is truly a combination of the stopping iteration and the choice of $\lambda$, so attention needs to be placed on finding a good stopping criterion.

\paragraph{Relation to existing methods} It is a well known result that, under $\bB=\bA^\mathsf{T}$, AB-/BA-GMRES are algebraically equivalent to LSQR and LSMR, respectively \cite{hansen2022gmres}. We now explore the corresponding relationship for hybrid AB-/BA-GMRES.  

First, we observe that hybrid LSMR, presented in \cite{chung2015hybrid}, and hybrid BA-GMRES are equivalent. By definition, they both find the following solution at iteration $k$:
\begin{equation*}
 \bx_k^{\text{h,A}^\mathsf{T}\text{A-GMRES}} = \arg\min_{\bx\in\mathcal{K}_k(\bA^\mathsf{T}\bA,\bA^\mathsf{T}\bb)}\| \bA^\mathsf{T}(\bA\bx-\bb)\|^2_2 + \lambda^2\| \bx\|^2_2 = \bx_k^{\text{h,LSMR}}.
\end{equation*}

Second, we show that hybrid LSQR and hybrid AB-GMRES are not equivalent. In particular,
at iteration $k$, the two methods produce different approximate solutions. Hybrid
AB-GMRES applied to the normal equations yields
\begin{align*}
\bx_k^{\text{h,}AA^\mathsf{T}\text{-GMRES}}
&= \bA^\mathsf{T}
\arg\min_{\bz \in \mathcal{K}_k(\bA\bA^\mathsf{T}, \bb)}
\left(
\| \bA\bA^\mathsf{T} \bz - \bb \|_2^2
+ \lambda^2 \| \bz \|_2^2
\right) \\
&= \arg\min_{\bx \in \mathcal{K}_k(\bA^\mathsf{T}\bA, \bA^\mathsf{T}\bb)}
\left(
\| \bA \bx - \bb \|_2^2
+ \lambda^2 \| (\bA^\mathsf{T})^{\dagger} \bx \|_2^2
\right),
\end{align*}
which, in general, differs from the hybrid LSQR iterate
$\bx_k^{\text{h,LSQR}}.$

\begin{figure}
\begin{minipage}{0.48\textwidth}
\begin{algorithm}[H] 
\caption{Hybrid AB-GMRES}\label{alg:bagmres1}
\begin{algorithmic}[1]
\State Choose $\bu_0$ and set $\bx_0=\bB\bu_0$
\State $\br_0 = \bb - \bA \bx_0$
\State $\bw_1 = \br_0 / \|\br_0\|_2$
\For{$k = 1, 2, \dots, K$}
    \State $\bq_k = \bA\bB \bw_k$
    \For{$i = 1, 2, \dots, k$}
        \State $\bh_{i,k} = \langle \bq_k, \bw_i \rangle$
        \State $\bq_k = \bq_k - \bh_{i,k} \bw_i$
    \EndFor
    \State $\bh_{k+1,k} = \|\bq_k\|_2$
    \State $\bw_{k+1} = \bq_k / \bh_{k+1,k}$
    \State Find optimal $\lambda$
    
    \State $\by_k \newline= \arg\min_{\by} \left\| \left(\smash{{}^{\;\;\bH_k}_{\lambda \bI}}\right)\by -\|\br_0\|_2\left(\smash{{}^{\be_1}_{0}}\right) \right\|_2$
    \State $\bx_k = \bx_0 + \bB\bW_k \by_k$
    \State $\br_k = \bb - \bA \bx_k$
    \State Stopping Rule
\EndFor
\end{algorithmic}
\end{algorithm}
\end{minipage}
\begin{minipage}{0.48\textwidth}
\begin{algorithm}[H]
\caption{Hybrid BA-GMRES}\label{alg:bagmres2}
\begin{algorithmic}[1]
\State Choose $\bx_0$
\State $\br_0 = \bB(\bb - \bA \bx_0)$
\State $\bw_1 = \br_0 / \|\br_0\|_2$
\For{$k = 1, 2, \dots, K$}
    \State $\bq_k = \bB\bA \bw_k$
    \For{$i = 1, 2, \dots, k$}
        \State $\bh_{i,k} = \langle \bq_k, \bw_i \rangle$
        \State $\bq_k = \bq_k - \bh_{i,k} \bw_i$
    \EndFor
    \State $\bh_{k+1,k} = \|\bq_k\|_2$
    \State $\bw_{k+1} = \bq_k / \bh_{k+1,k}$
    \State Find optimal $\lambda$
    \State $\by_k \newline= \arg \min_{\by} \left\| \left(\smash{{}^{\;\;\bH_k}_{\lambda \bI}}\right)\by -\|\br_0\|_2\left(\smash{{}^{\be_1}_{0}}\right) \right\|_2$
    \State $\bx_k = \bx_0 + \bW_k \by_k$
    \State $\br_k = \bB(\bb - \bA \bx_k)$
    \State Stopping Rule
\EndFor
\end{algorithmic}
\end{algorithm}
\end{minipage}
\end{figure}
\subsection{Regularization Parameter Selection}
\label{sec:regparam}
An important and often challenging aspect of regularization methods is the selection of an appropriate regularization parameter $\lambda$. The quality of the reconstructed solution is highly sensitive to this choice: under-regularization may lead to noise amplification and instability, while over-regularization can oversmooth the solution and suppress important structural features \cite{hansen1998rank, hansen2006deblurring}. 
In this work, we investigate two widely used and well-established parameter-selection techniques: the L-curve criterion \cite{Hansenlc} and generalized cross validation (GCV) \cite{golub1997generalized}. Both methods aim to balance data fidelity and solution regularity, but they differ in their underlying principles and practical behavior.

In the context of hybrid Krylov subspace methods, the regularization parameter $\lambda$ be computed efficiently for each of the small projected problems \eqref{eq:hybrid_proj} at each iteration. This allows the regularization parameter to be selected adaptively as the Krylov subspace grows. For this reason, we maintain the notation of the problem in \eqref{eq:hybrid_proj}.

\paragraph{L-curve Criterion}
The L-curve criterion is a graphical method that characterizes the trade-off between the residual norm and the solution norm in Tikhonov regularization \cite{hansen1992analysis,hansen2010discrete}. Specifically, the L-curve for the projected problem at iteration $k$ is defined as the parametric plot of
\begin{equation}
\left( \log \|\bH_k\by_\lambda - \|\br_0\|_2 \be_1\|_2, \; \log \|\by_\lambda\|_2 \right),
\end{equation}
where $\bx_\lambda$ denotes the Tikhonov-regularized solution corresponding to the parameter $\lambda$. The curve typically exhibits an ``L''-shaped structure, with the corner representing a compromise between minimizing the residual norm and controlling the solution norm. While the L-curve criterion is intuitive and often effective, its performance can be sensitive to noise level and problem conditioning, and identifying the corner of the curve may be challenging in practice for severely ill-posed problems.

\paragraph{Generalized Cross Validation}
Generalized cross validation (GCV) provides an automated and statistically motivated approach for selecting the regularization parameter \cite{wahba1990spline,hansen2010discrete}. The GCV criterion seeks to minimize a predictive error estimate that approximates leave-one-out cross validation, without requiring explicit knowledge of the noise level. For Tikhonov regularization in the projected problem at iteration $k$, the GCV function is given by
\begin{equation}
\mathrm{GCV}(\lambda) =
\frac{\|\bH_k\by_\lambda - \|\br_0\|_2 \be_1\|_2^2}
{\left( \mathrm{trace}(\bI - \bH_k \bH_{k,\lambda} \right)^2},
\end{equation}
where $\bH_{k,\lambda}= (\bH_k\t \bH_k +\lambda^2 \bI)^{-1} \bH_k\t$ denotes the regularized projected inverse operator.

The GCV criterion is fully automatic and does not rely on an heuristic interpretation, which makes it attractive for practical implementations. However, in the presence of correlated noise or strong model mismatch, the GCV criterion may underestimate the optimal regularization parameter, leading to mild under-regularization. The GCV function from TRIPs-Py is used in this paper \cite{pasha2024trips}.

\paragraph{Discussion}
Both the L-curve and GCV methods are well suited for hybrid ABBA-GMRES, as they operate on the small projected problems generated during the Arnoldi process. This enables efficient, iteration-dependent parameter selection without additional forward or backprojection operations. In our numerical experiments, we observe that both criteria successfully mitigate semi-convergence and improve reconstruction stability compared to non-hybrid methods. While the L-curve criterion provides valuable insight into the trade-off between data fidelity and regularization, GCV offers a more automated and scalable alternative. The choice between these methods ultimately depends on the noise characteristics, problem conditioning, and computational constraints. 
\subsection{Stopping Criteria}\label{sec:stop_critieria}

Determining when to terminate an iterative method is a critical component of solving inverse problems. Stopping too early may result in an under-resolved solution that fails to capture essential features, while stopping too late can amplify noise and numerical errors. In the context of hybrid methods, the risk of over-iteration is mitigated by the explicit regularization term in the projected problems.
In this work, we consider three stopping criteria: the Discrepancy Principle (DP), the Normalized Cumulative Periodogram (NCP), and Residual Norm Stagnation (RNS).

\paragraph{Discrepancy Principle}

The DP seeks to terminate the iterations once the residual is consistent with the level of noise in the data. Consider the perturbed system
\[
\bA \bx = \bb + \be,
\]
where $\be$ represents the measurement noise. The DP prescribes stopping at iteration $k$ when the residual norm satisfies is approximately the norm of the noise $\|\be\|_2$. 
In practice, a safety factor $\tau > 1$ is introduced, and the stopping criterion becomes
\[
\|\bA \bx_k - \bb\|_2 \geq \tau \|\be\|_2.
\]
A key limitation of the DP is that it requires an estimate of the noise level $\|\be\|_2$, which is often unavailable or difficult to obtain accurately.

\paragraph{Normalized Cumulative Periodogram}

The NCP is conceptually related to the DP in that it aims to detect the point at which the residual consists predominantly of noise. Unlike the DP, however, the NCP does not require prior knowledge of the noise norm \cite{hansen2010discrete}. Instead, it analyzes the frequency content of the residual by examining its power spectrum and identifies the iteration at which the residual behaves like white noise. While effective, this approach is computationally more involved than norm-based criteria.

\paragraph{Residual Norm Stagnation}
The RNS is based on monitoring the relative change in the residual norm between successive iterations. The iterations are terminated when
$$ \frac{\left|\|\bA \bx_{k-1} - \bb\|_2 - \|\bA \bx_k - \bb\|_2\right|}{\|\bA \bx_{k-1} - \bb\|_2} < \varepsilon,
$$ where $\varepsilon$ is a prescribed tolerance. This criterion assumes that once the residual norm ceases to decrease significantly, further iterations are unlikely to yield meaningful improvements in the solution.

\subsection{Restarting} As Krylov subspace methods progress, the dimension of the underlying subspace increases with each iteration, leading to growing memory requirements and higher computational costs due to orthogonalization. To control these costs, it is common to employ \emph{restarting}, in which the iterative process is terminated after a fixed number of iterations and then re-initiated, see for instance \cite{saad2003iterative, greenbaum1997iterative, buccini2023limited, simoncini2000convergence, morgan2002gmres}.

Specifically, after $p$ iterations, the algorithm is restarted using as an initial guess the approximate solution obtained at the end of the previous cycle. This effectively limits the dimension of the Krylov subspace to at most $p$, thereby reducing memory usage and improving computational efficiency. Restarting is particularly important for large-scale problems where storage and computational constraints are significant \cite{saad1986gmres}.

While restarting helps manage resources, it may also slow convergence or lead to the loss of important spectral information accumulated in the Krylov subspace \cite{greenbaum1997iterative}. Hybrid methods partially mitigate this drawback by incorporating regularization within each restart cycle, allowing meaningful solution updates even when the subspace is truncated \cite{kilmer2006hybrid, chung2008hybrid}. In practice, the choice of the restart parameter $p$ represents a trade-off between computational efficiency and convergence behavior, and is typically selected based on problem size, noise level, and available computational resources. For very large scale or massive data, this is the only computationally feasible way of implementing AB-/BA-GMRES or their hybrid counterparts. However, ignoring the previous solution space leads to discarding important information, necessitating additional work to reintroduce the information into the solution space. Recycling is known as a more effective method that alternately enlarges and judiciously compresses the solution subspace, keeping the dimension bounded independent of the number of iterations, while still preserving the relevant information throughout the process for fast convergence \cite{pasha2023recycling, jiang2021hybrid}. In this paper we only consider restarting and leave recycling as a potential direction to be explored in the future.

\section{Computational Results}
\label{sec:computational}
In this section, we investigate the numerical performance of the proposed hybrid ABBA-GMRES methods through a series of representative computational experiments. The primary objective of these experiments is to assess the effectiveness of hybrid regularization in stabilizing convergence and improving reconstruction quality when unmatched forward and backprojector pairs are employed.

First, we consider a small synthetic computed tomography example which has been previously used specifically to evaluate the performance of iterative solvers for unmatched projector pairs \cite{website_unmatched}. This toy example is designed to compare different methods in both simulated unmatched and matched projector  scenarios. Its primary purpose is to examine the behavior of the proposed hybrid AB-/BA-GMRES methods in comparison with hybrid LSQR and hybrid LSMR, and to illustrate the theoretical results presented in the paper. To isolate the effects of the algorithms themselves, we fix the regularization parameter throughout the experiment. Since the native language of the example in \cite{website_unmatched} is MATLAB, as well as the hybrid LSMR codes distributed by their authors \cite{website_hybrid_lsmr}, we use this programming language for this particular experiment. The codes used for LSQR and hybrid LSQR, are obtained using IRtools \cite{IRtools}. We emphasize that this small test problem is only designed to validate theoretical findings and not to reflect efficiency or performance.

Second, we consider two test problems that are designed to reflect key challenges encountered in large-scale X-ray computed tomography, including more realistic problems and the need of efficient GPU-specific implementations. In particular, the data is generated using the ASTRA Toolbox \cite{ASTRA} using GPU acceleration. In this case, Python is employed to generate larger-scale, more realistic examples.

For each problem, we compare the behavior of the hybrid AB-/BA-GMRES methods against their non-hybrid counterparts, with particular emphasis on convergence characteristics, sensitivity to regularization parameter selection, and robustness with respect to semi-convergence. Specifically, we examine the influence of different regularization parameter selection strategies, including the L-curve criterion and generalized cross validation, and evaluate how these choices affect convergence behavior and reconstruction quality.

In addition to reconstruction accuracy, we assess the computational efficiency of the proposed methods. Since forward and backprojection operations dominate the overall cost in large-scale CT, we report iteration counts, timing results, and convergence rates to demonstrate that hybrid ABBA-GMRES methods achieve improved stability with minimal additional computational overhead. Together, these experiments provide a comprehensive evaluation of the proposed methods and illustrate their practical advantages for large-scale tomographic reconstruction with unmatched backprojectors.
Numerical experiments were run on a personal computer with specifications summarized in Table \ref{tab:comp_spec}.
In our experiments, we assess the quality of the reconstructions using the relative reconstruction error (RRE) between the recovered $k$-th iterate $\bx_{k}$ and the ground truth image $\bx_{\text{true}}$ defined as $\text{RRE}(\bx_{k}, \bx_{\text{true}}) = \frac{\|\bx_{k} - \bx_{\text{true}}\|_2}{\|\bx_{\text{true}}\|_2}$ 
and structural similarity index (SSIM) measures \cite{wang2004image}. 
All test problems and algorithm implementations in Python will be made publicly available at 
\begin{center}
\href{https://github.com/rbentley5/GPU-CT-Unmatched-back-projectors}{https://github.com/rbentley5/Hybrid\_ABBA\_GMRES} 
\end{center}
once the manuscript is accepted to the journal. Other implementations of ABBA methods developed in C++ are incorporated into the Re::Solve software package \cite{swirydowicz2023resolve}. The implementations in Re::Solve focus on GPU implementations of ABBA GMRES to be used for problems involving sparse matrices.

\begin{table}[h]
    \centering
    \vspace{0.5em}
    \begin{tabular}{ll}
        \hline
        Memory: & 32 GB, 4800 MHz, DDR4 \\
        CPU: & 12th Gen Intel(R) Core(TM) i7-12800H, 14 cores \\
        GPU Memory: & 4 GB dedicated, 16 GB shared \\
        GPU: & NVIDIA RTX A2000 Laptop GPU \\
        \hline
    \end{tabular}
    \caption{Computer specifications.}
    \label{tab:comp_spec}
\end{table}
\subsection{Test Problem 1}
This toy problem is introduced solely to assess selected theoretical properties of the proposed method, and it is taken from \cite{website_unmatched}. In particular, this simulates a small CT example, where $\bA$ is a sparse matrix –with $98.6\%$ sparsity– that represents the forward projector corresponding to 180 detector angles and 128 detector pixels but the original source does not specify the geometry further. The solution image is 128 $\times$ 128 pixels, and additive white Gaussian noise with a low noise level or $0.1\%$ has been added to the measurements. In this example, we use MATLAB implementations to enable comparison with other methods
that are likewise provided in MATLAB \cite{IRtools,website_hybrid_lsmr}.
The problem size is sufficiently small that the explicit transpose of the system matrix can be computed and used as a reference. In this setting, we investigate the behavior of AB-/BA-GMRES and its hybrid variants, and compare their performance with that of LSQR and LSMR, as well as their corresponding hybrid counterparts. This comparison allows us to isolate the differences between the methods, as well as their equivalences, in a controlled environment, and it is not meant to be a realistic example.

\begin{figure}[ht]
  \centering
  \begin{minipage}[t]{0.49\textwidth}
    \centering
    \includegraphics[height=5.5cm]{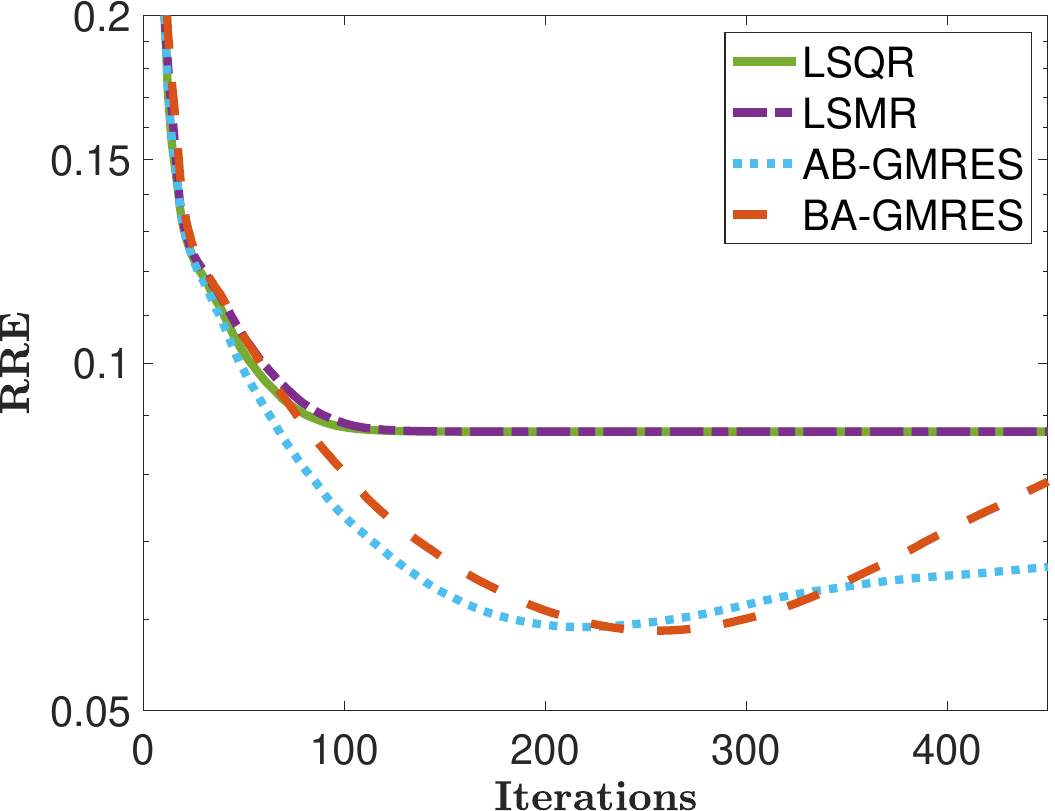}
  \end{minipage}%
  \begin{minipage}[t]{0.49\textwidth}
    \centering
    \includegraphics[height=5.5cm]{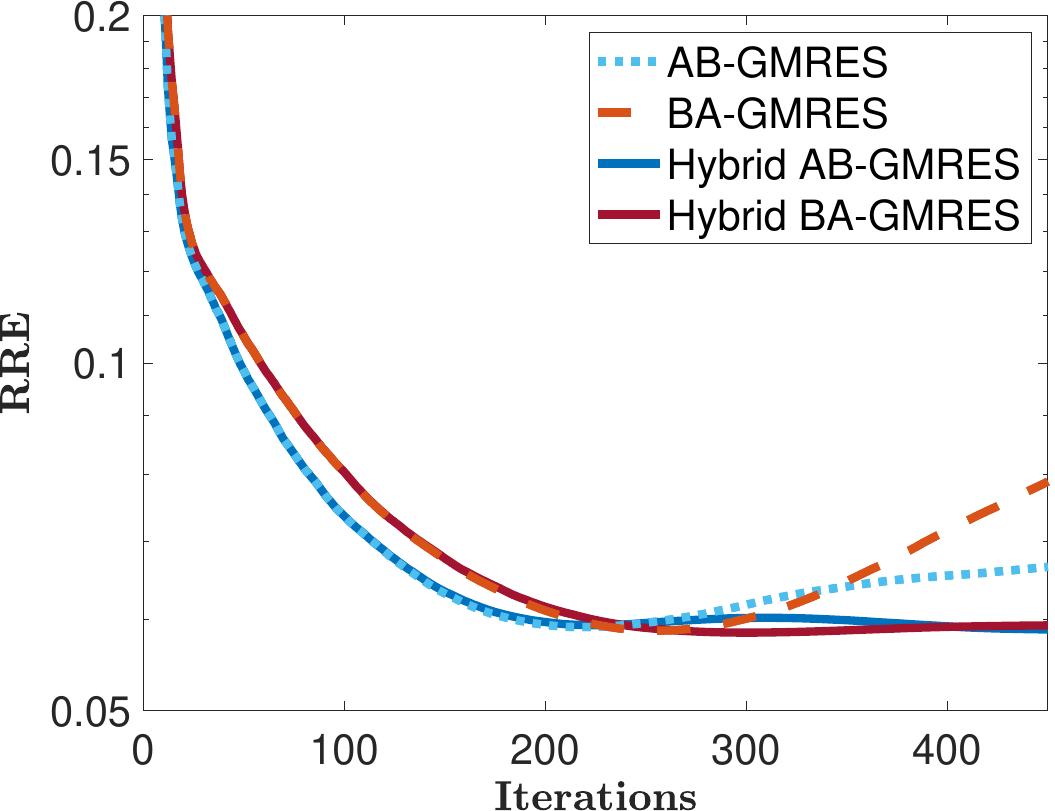}
  \end{minipage}
  \caption{Relative reconstruction error norms for different solvers using an unmatched backprojector, i.e. $\bB\neq\bA\t$.} \label{Fig:Ex1_1}
\end{figure}

In Figure \ref{Fig:Ex1_1} (left panel), we can observe the difference in the performance of LSQR and LSMR for this unmatched example, where the implementations of LSQR and LSMR are taken from \cite{IRtools,website_hybrid_lsmr}, repectively, and run with re-orthogonalization. One can observe that the AB-/BA-GMRES methods perform much better for this example. However, both AB-/BA-GMRES display semi-convergence, while LSQR and LSMR do not for this example (this might be related to stagnation generated by the unmatched backprojector, but it is not in the scope of this paper to investigate this further). In Figure \ref{Fig:Ex1_1} (right panel), the  AB-/BA-GMRES methods are displayed along with their hybrid counterparts, with regularization parameters chosen manually to stabilize convergence ($\lambda=0.0015$ and $\lambda=0.005$, respectively). In the next experiments, we show automatic ways of selecting this.

\begin{figure}[H]
  \centering
  \begin{minipage}[t]{0.49\textwidth}
    \centering
    \includegraphics[height=5.5cm]{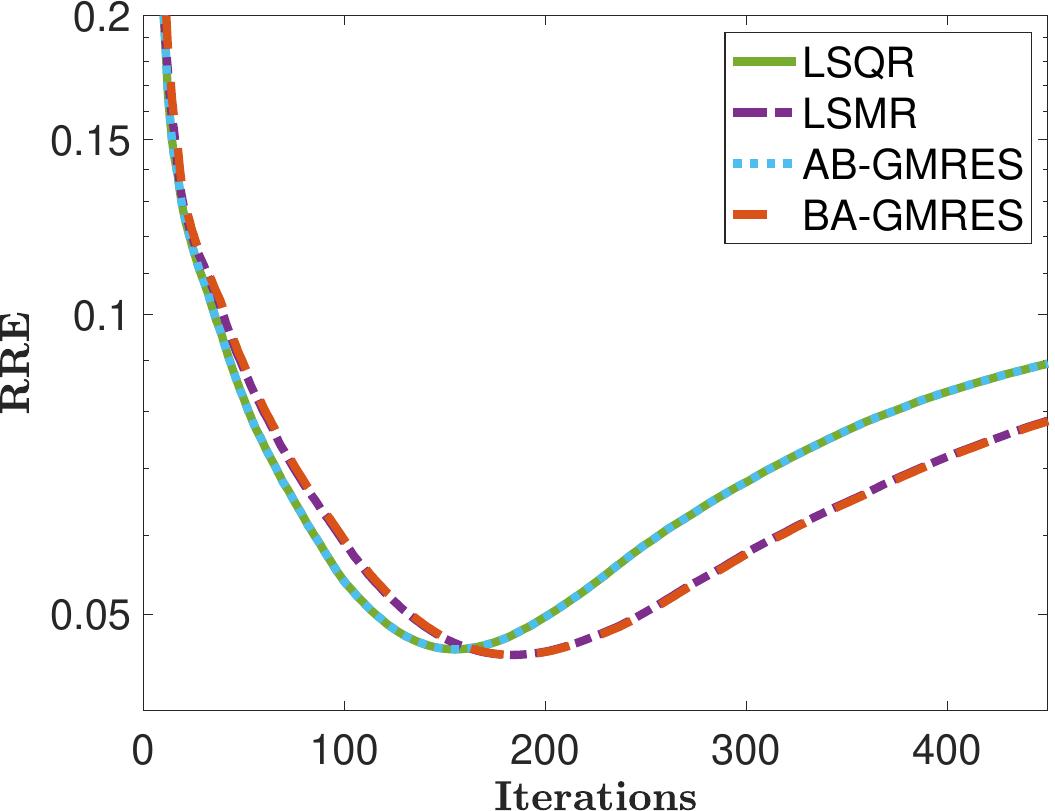}
  \end{minipage}%
  \begin{minipage}[t]{0.49\textwidth}
    \centering
    \includegraphics[height=5.5cm]{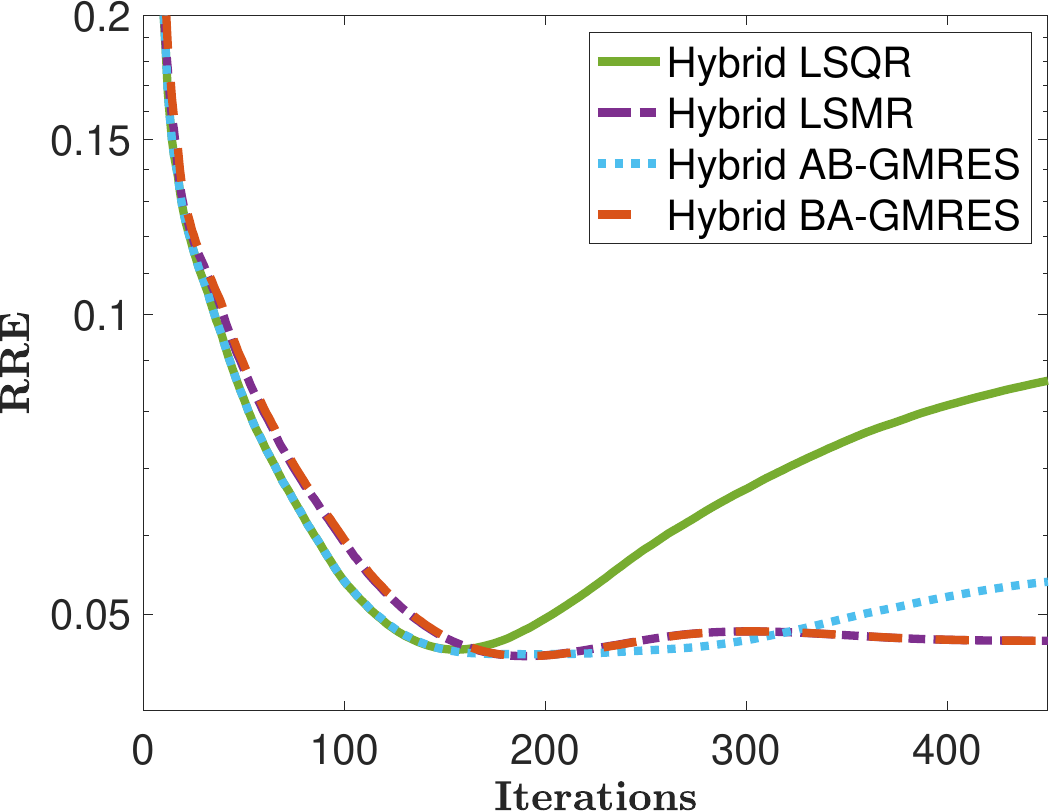}
  \end{minipage}
\caption{Relative reconstruction error norms for different solvers using the same system matrix $\bA$ and a matched backprojector, i.e. $\bB=\bA\t$.}
\label{Fig:Ex1_2}
\end{figure}

Figure \ref{Fig:Ex1_2} is displayed to exemplify the theory derived for the matched case, i.e. where we impose that $\bB=\bA\t$. As one can expect, LSQR and AB-GMRES are equivalent in this case, as well as LSMR and BA-GMRES. However, this is only true for LSMR in the hybrid case. 

\subsection{Test Problem 2}
\begin{figure}[h]
  \centering
  \begin{minipage}[t]{0.44\textwidth}
    \centering
    \includegraphics[height=5.5cm]{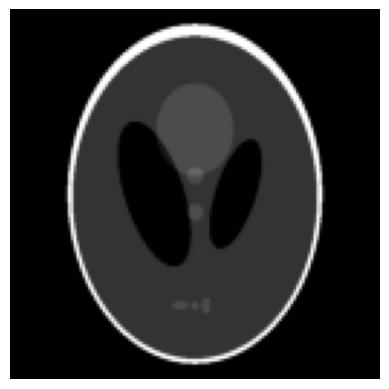}
  \end{minipage}%
  \begin{minipage}[t]{0.5\textwidth}
    \centering
    \includegraphics[height=3cm]{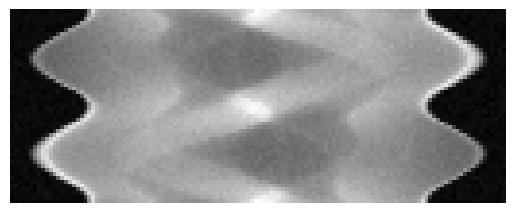}
  \end{minipage}
    \caption{Left panel: True image of size $128\times 128$. Right panel: noisy sinogram from 50 view angles and 2.5\% Gaussian noise.}
    \label{Fig_Ex2_settting}
\end{figure}

The second test problem involves a 128 x 128 Shepp-Logan phantom, in a synthetic setting that is more realistic than the previous example. The sinogram was obtained using parallel geometry through 50 view angles with an noise level of 2.5\% yielding $\bb\in \R^{6400}$. The forward projector, $\bA \in \R^{6400\times16384}$, and back projector, $\bB \in \R^{16384\times6400}$, are generated using ASTRA Toolbox \cite{ASTRA}. Both the image and the corrupted measurements can be observed in Figure \ref{Fig_Ex2_settting}. For this example, we do an exhaustive comparison, testing the effect of using different regularization parameter selection and stopping criteria, restarting techniques, and quantify efficiency.

\begin{figure}[h]
  \centering
  \begin{minipage}[t]{0.5\textwidth}
    \centering
    \includegraphics[height=5.5cm]{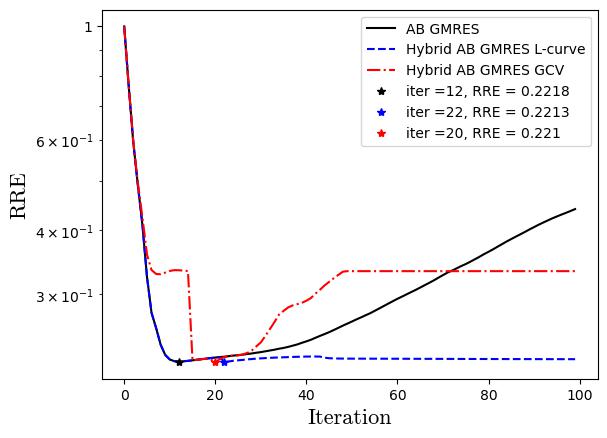}
  \end{minipage}%
  \begin{minipage}[t]{0.5\textwidth}
    \centering
    \includegraphics[height=5.5cm]{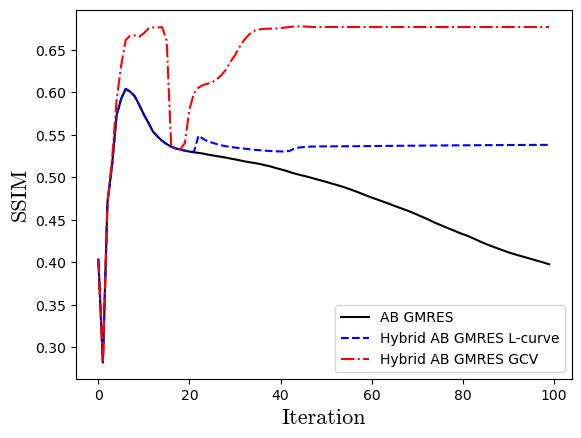}
  \end{minipage}
  \caption{Relative reconstruction error norms and SSIM for AB-GMRES with different regularization parameter choices. The markers correspond to the iterations with minimum RRE for each method.}\label{Fig:Reg_param_AB}
\end{figure}

\begin{figure}[h]
  \centering
  \begin{minipage}[t]{0.5\textwidth}
    \centering
    \includegraphics[height=5.5cm]{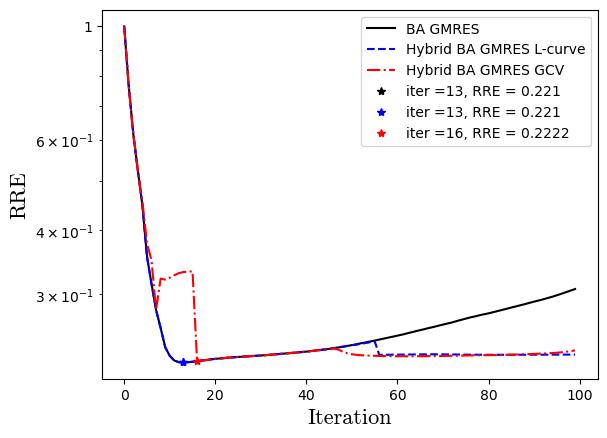}
  \end{minipage}%
  \begin{minipage}[t]{0.5\textwidth}
    \centering
    \includegraphics[height=5.5cm]{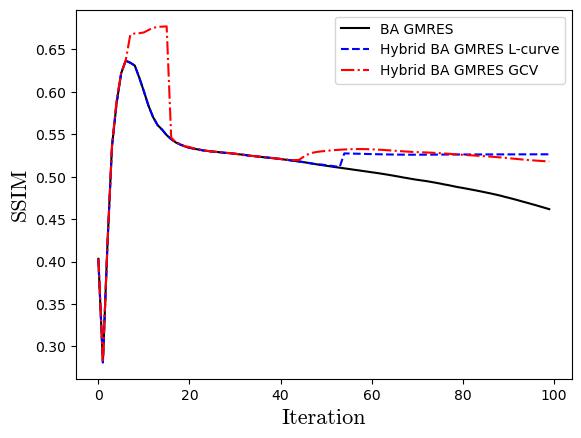}
  \end{minipage}
    \caption{Relative reconstruction error norms and SSIM for BA-GMRES with different regularization parameter choices. The markers correspond to the iterations with minimum RRE for each method.}\label{Fig:Reg_param_BA}
\end{figure}

\paragraph{Regularization Parameter Selection}
Figures \ref{Fig:Reg_param_AB} and  \ref{Fig:Reg_param_BA} show the relative reconstruction error (RRE) norms and SSIM of hybrid AB-GMRES and hybrid BA-GMRES, respectively, against the number of iterations, using different regularization parameter criteria. 

In Figure \ref{Fig:Reg_param_AB} (left panel), one can observe that the relative error norm using hybrid AB-GMRES completely stabilizes at the minimum error when using the L-curve, while GCV semi-converges a bit before stabilizing. However, in Figure \ref{Fig:Reg_param_AB} (right panel), the SSIM using hybrid AB-GMRES stabilizes at the maximum for GCV while L-curve stabilizes at a lower value. In Figure \ref{Fig:Reg_param_BA} (left panel), we observe that both regularization parameter choices are very effective in terms of the RRE in this example for BA-GMRES. In Figure \ref{Fig:Reg_param_BA} (right panel) both L-curve and GCV show very similar behavior, both attaining a maximum before stabilizing at a lower SSIM. Using the GCV criterion we obtain a higher maximum SSIM, similar to hybrid AB-/BA-GMRES methods. It is worth noting that the original BA-GMRES shows less severe semi-convergence than AB-GMRES for this example.  

\begin{figure}[h]
    \centering
    \includegraphics[width=0.5\linewidth]{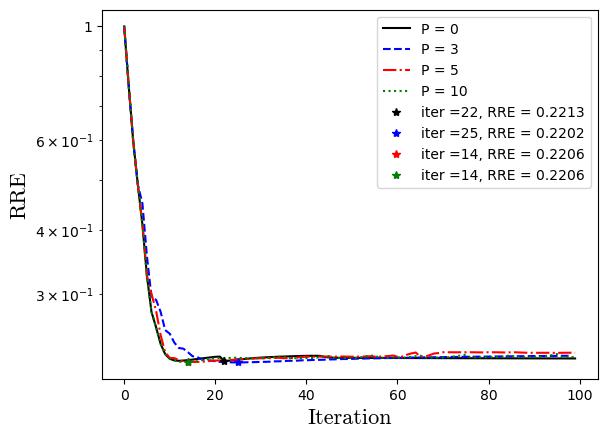}
    \caption{Relative reconstruction error norm against the number of iterations for Hybrid AB-GMRES with L-curve restarted every $p$ iterations. The markers correspond to the iterations with minimum RRE for each method.}\label{Fig:restarts}
\end{figure}

\paragraph{Restarted Iterations}
Figure \ref{Fig:restarts} shows the relative error norm of the solution obtained using AB-GMRES restarted after a different amount of iterations. As expected from the theory, for this problem, decreasing $p$ leads to a slower convergence in terms of number of iterations. However, there is possibility for a speedup in computational time depending on the problem size. 
\paragraph{Resource Usage}
Table \ref{tab:usage2} summarizes different resources usage metrics for this example using different solvers. One can observe that the AB-GMRES methods unilaterally performed better in terms of time and memory usage. The L-curve is more efficient than GCV in all cases especially for the non-restarted methods. Restarted methods take less memory and are faster in general. In terms of GPU usage, there does not seem to be any clear trends.

\begin{table}[H]
    \centering
    \begin{tabular}{|c|c|c|c|c|c|}
        \hline
      Method   & Regularization Param. &  Restarted It. & Time (sec) & RAM (mb) & GPU Usage\\
      AB & None & 0 & 0.53 & 18.9 & 0.053\\
      AB & None & 10 & 0.36 & 10.4 & 0.093\\
      AB & L-curve & 0 & 0.82 & 18.9 & 0.07\\
      AB & L-curve & 10 & 0.54 & 10.4& 0.06\\
      AB & GCV & 0 & 4.27 & 19.2 & 0.04\\
      AB & GCV & 10 & 1.43 & 10.5 & 0.043\\
\hline
      BA & None & 0 & 0.85 & 23.0& 0.06\\
      BA & None & 10 & 0.36 &11.4 & 0.04\\
      BA & L-curve & 0 & 1.09 & 23.0& 0.06\\
      BA & L-curve & 10 & 0.66 & 11.4& 0.077\\
      BA & GCV & 0 & 4.88 & 23.4 & 0.027\\
      BA & GCV & 10 & 1.64 & 11.5 & 0.033\\
\hline
      
    \end{tabular}
    \caption{Resource Usage for Test Problem 2}
    \label{tab:usage2}
\end{table}
\subsection{Test Problem 3}
\begin{figure}[h]
  \centering
  \begin{minipage}[t]{0.47\textwidth}
    \centering
    \includegraphics[height=6.5cm]{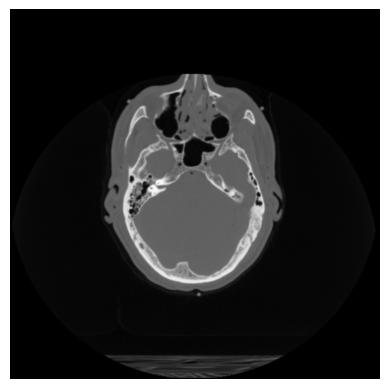}
  \end{minipage}%
  \begin{minipage}[t]{0.53\textwidth}
    \centering
    \includegraphics[height=.75cm]{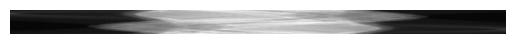}
  \end{minipage}
    \caption{Left panel: True image of size $1024\times 1024$ pixels. Right panel: noisy sinogram from 50 view angles and 1.5\% Gaussian noise.}
    \label{Fig:Ex3_setting}
\end{figure}

The third test problem corresponds to a 1024 x 1024 Logan phantom. The sinogram was obtained using fan-flat geometry through 50 view angles with a noise level of 1.5\% yielding $\bb\in \R^{51200}$. The forward projector, $\bA \in \R^{51200\times1048576}$, and back projector, $\bB \in \R^{1048576\times51200}$, are generated using the ASTRA Toolbox \cite{ASTRA}. We show the exact solution as well as the corrupted measurements in Figure \ref{Fig:Ex3_setting}. In this example, we compare the performance of different regularization parameter selection methods as well as stopping criteria, and we evaluate quantitatively the efficiency of the new methods using different metrics.

\begin{figure}[h]
  \centering
  \begin{minipage}[t]{0.49\textwidth}
    \centering
    \includegraphics[height=5.5cm]{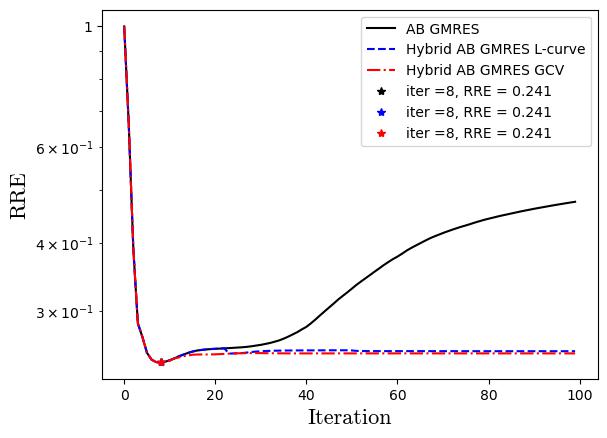}
  \end{minipage}%
  \begin{minipage}[t]{0.49\textwidth}
    \centering
    \includegraphics[height=5.5cm]{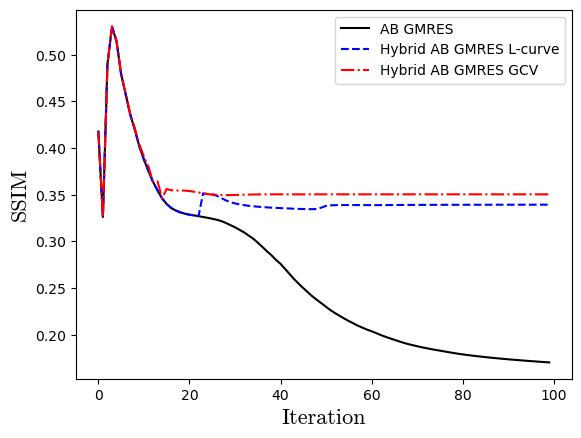}
  \end{minipage}
  \caption{Relative error norms and SSIM for AB-GMRES with different regularization parameter choices. The markers correspond to the iterations with minimum RRE for each method.}\label{Fig:reg_param_ABGMRES_2}
\end{figure}

\begin{figure}[H]
  \centering
  \begin{minipage}[t]{0.49\textwidth}
    \centering
    \includegraphics[height=5.5cm]{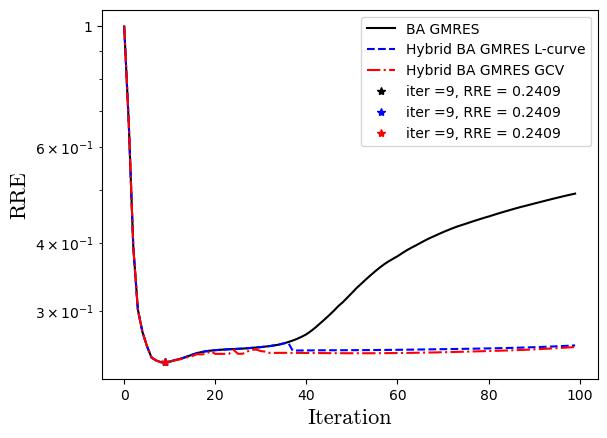}
  \end{minipage}%
  \begin{minipage}[t]{0.49\textwidth}
    \centering
    \includegraphics[height=5.5cm]{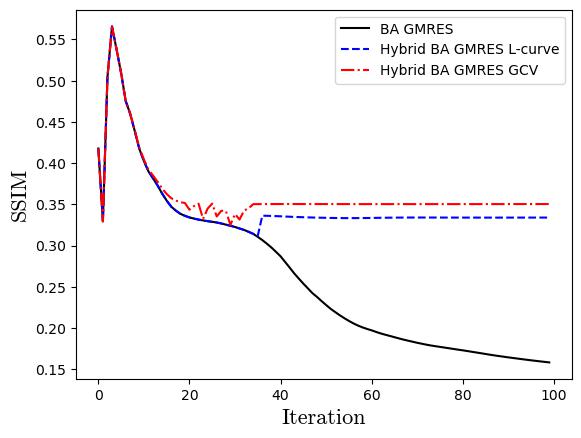}
  \end{minipage}
  \caption{Relative reconstruction error  norms and SSIM for BA-GMRES with different regularization parameter choices. The markers correspond to the iterations with minimum RRE for each method.}\label{Fig:reg_param_BAGMRES_2}
\end{figure}

\paragraph{Regularization Parameter Selection}
Figures \ref{Fig:reg_param_ABGMRES_2} and  \ref{Fig:reg_param_BAGMRES_2} show the relative error norms and SSIM of hybrid AB-GMRES and hybrid BA-GMRES, respectively, against the number of iterations, using different regularization parameter criteria. 

In Figure \ref{Fig:reg_param_ABGMRES_2} (left panel), we observe that the relative error norm obtained with hybrid AB-GMRES behaves very similarly when using either the GCV or the L-curve criterion, exhibiting only a slight degree of semi-convergence before stabilizing. This is also supported by the SSIM in Figure \ref{Fig:reg_param_ABGMRES_2}: both regularization parameter selection methods attain the same maximum SSIM, after which the value decreases slightly before stabilizing at a similar level. In Figure \ref{Fig:reg_param_BAGMRES_2} (left panel), we observe that both regularization parameter choices exhibit slightly less stability in terms of the RRE than AB-GMRES, as evidenced by a slight increase around 90 iterations. In Figure \ref{Fig:reg_param_BAGMRES_2} (right panel) both regularization parameter selection methods show a very similar behavior, attaining a maximum before stabilizing at a lower SSIM.

\begin{figure}[h]
  \centering
  \begin{minipage}[t]{0.49\textwidth}
    \centering
    \includegraphics[height=5.5cm]{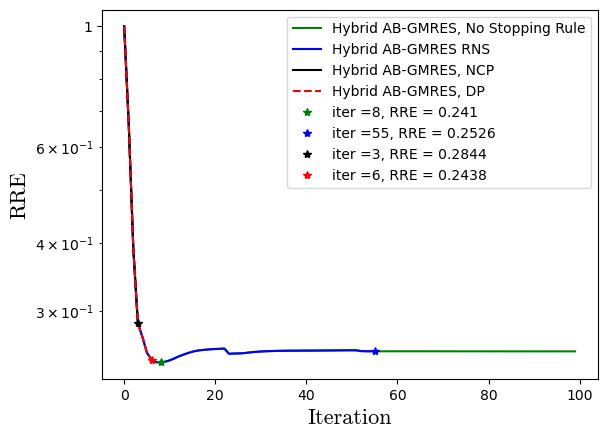}\\
  \end{minipage}
  \begin{minipage}[t]{0.49\textwidth}
    \centering
    \includegraphics[height=5.5cm]{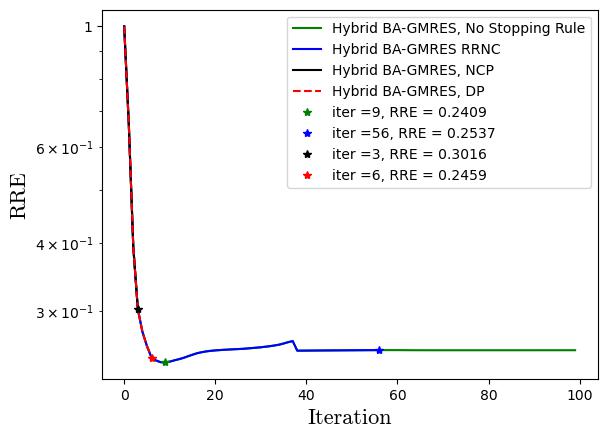}
  \end{minipage}
  \caption{Relative reconstruction errors for hybrid AB-/BA-GMRES using the L-curve regularization parameter choice criterion, and different stopping criteria. The markers correspond to the stopping iteration for each of the criteria, where the green marker corresponds to the optimal stopping iteration (i.e. the one minimizing the RRE).}
  \label{Fig:ABBA_stopping}
\end{figure}

\begin{figure}[H]
    \centering
    \includegraphics[width=0.9\linewidth]{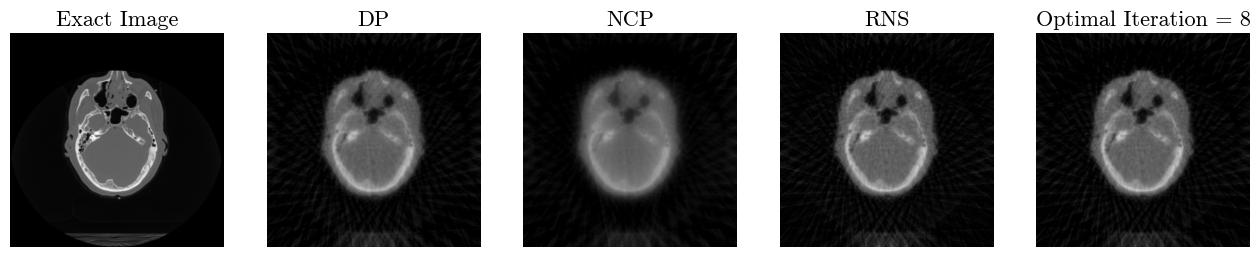}
    \caption{Reconstructed images using hybrid AB-GMRES with different stopping criterion.}
    \label{fig:rec_ex31}
\end{figure}

\begin{figure}[h]
    \centering
    \includegraphics[width=0.9\linewidth]{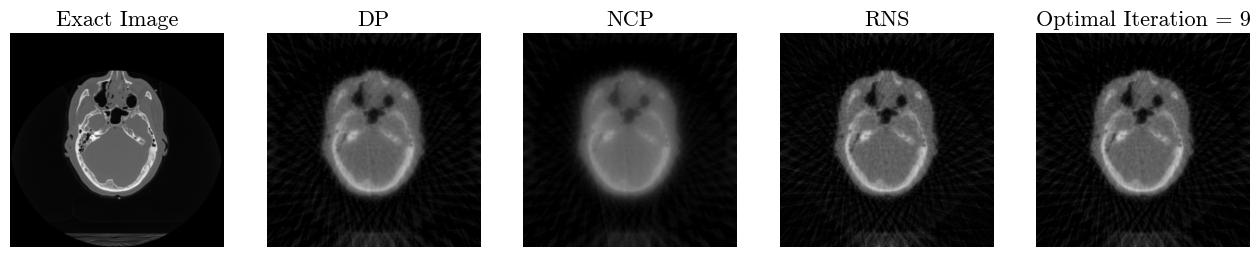}
    \caption{Reconstructed images using hybrid BA-GMRES with different stopping criterion.}
    \label{fig:rec_ex32}
\end{figure}

\paragraph{Stopping criteria} 
Figure \ref{Fig:ABBA_stopping} shows the relative reconstruction errors for hybrid AB-/BA-GMRES using a fixed regularization parameter selection criterion (the L-curve) and different stopping criteria. In this example, we observe that the DP stopping criterion performs very well for both hybrid AB-GMRES and BA-GMRES. On the other hand, the NCP method stops too early for both methods, while the RNS criterion stops too late. It is important to remark that the RNS is expected to have a delayed stopping iteration; however, given the nature of hybrid methods, this only decreases performance (by requiring more iterations) and does not sacrifice reconstruction quality (due to the effect of the explicit regularization). The reconstruction using AB-GRMES and BA-GMRES given the different stopping criteria can be observed in Figure \ref{fig:rec_ex31} and \ref{fig:rec_ex32}.

\paragraph{Resource Usage} Table \ref{tab:ex2} shows different performance metrics for this example. In this case, the AB-GMRES methods are generally much more efficient, especially the non-restarted version, compared to the BA-GMRES methods. This is due to the problem being a highly undetermined system. 
For BA-GMRES, this leads to much larger matrices compared to AB-GMRES. We can also observe that using the L-curve as a regularization parameter choice is only a bit more computationally demanding than not having any regularization, and it is much more efficient than using GCV. As expected, we can observe that restarting leads to faster times and less memory usage.

\begin{table}[h]
    \centering
    \begin{tabular}{|c|c|c|c|c|c|}
        \hline
      Method   & Regularization Param. &  Restarted It. & Time (sec) & RAM (mb) & GPU Usage\\
      AB & None & 0 & 6.03 & 897.6 & 0.39\\
      AB & None & 10 & 4.04 & 501.5 & 0.46\\
      AB & L-curve & 0 & 7.07 & 897.6 & 0.36\\
      AB & L-curve & 10 & 4.25 & 501.5 & 0.45\\
      AB & GCV & 0 & 14.56 & 897.8 & 0.42\\
      AB & GCV & 10 & 5.59 & 501.6 & 0.40\\
\hline
      BA & None & 0 & 47.61 & 1308.7 & 0.41\\
      BA & None & 10 & 7.05 & 586.9 & 0.29\\
      BA & L-curve & 0 & 49.31 & 1308.7 & 0.40\\
      BA & L-curve & 10 & 6.90 & 586.9 & 0.30\\
      BA & GCV & 0 & 53.53 & 1308.9 & 0.41\\
      BA & GCV & 10 & 8.47 & 587 & 0.29\\
\hline
\end{tabular}
\caption{Resource Usage for Test Problem 3}
\label{tab:ex2}
\end{table}

\section{Conclusions}

In this work, we proposed hybrid versions of AB-GMRES and BA-GMRES for large-scale X-ray computed tomography with unmatched forward and backprojector pairs. By incorporating Tikhonov regularization directly into the Krylov subspace iterations, the proposed hybrid AB-/BA-GMRES framework explicitly addresses the semi-convergence behavior that limits the effectiveness of standard AB-/BA-GMRES methods. Moreover, we extended existing theory relating the AB-/BA-GMRES methods to LSQR and LSMR under matched conditions, showing that hybrid BA-GMRES is equivalent to hybrid LSMR, while hybrid AB-GMRES is not equivalent to hybrid LSQR. This is due to the fact that these are `project-then-regularize' methods, and in this context, this s not equivalent to `regularize-then-project'.

Through a series of computational experiments, we demonstrated  that using hybrid AB-/BA-GMRES methods significantly improves convergence stability and reconstruction quality compared to their standard counterparts. In particular, the introduction of explicit regularization mitigates noise amplification at later iterations and reduces sensitivity to heuristic stopping rules. Among the regularization parameter selection strategies considered, the L-curve criterion consistently produced the most efficient reconstructions, while generalized cross validation offered a fully automatic alternative with a competitive performance.

We also showed that restarting the hybrid methods every $p$ iterations effectively controls Krylov subspace growth and limits memory requirements, leading to improved computational efficiency without sacrificing reconstruction quality. Although the inclusion of Tikhonov regularization introduces an additional overhead due to having to solve regularized projected problems, this cost is negligible relative to the dominant expense of forward and backprojection operations. Overall, the hybrid methods achieved improved robustness at the cost of only moderate increases in runtime.

Finally, we observed that stopping criteria based on the relative change of the residual provided more consistent and reliable termination behavior than classical discrepancy-based or norm-correction principles, particularly in the presence of unmatched operators and unknown noise levels. This further underscores the practical advantages of the hybrid AB-/BA-GMRES framework for real-world CT applications.

Taken together, these results indicate that hybrid AB-/BA-GMRES offers a stable, scalable, and effective approach for tomographic reconstruction with unmatched backprojectors. Future work will focus on extending the framework to incorporate more general regularization operators, adaptive restarting strategies, and applications to three-dimensional and time-resolved imaging problems.

\section*{Acknowledgments}
MP, RB and JZH acknowledge support from NSF DMS 2410699.  Any opinions, findings, conclusions, or recommendations expressed in this material are those of the authors and do not necessarily reflect the views of the National Science Foundation. RB and JZH acknowledge support from the Computational Modeling and Data Analytics undergraduate research funds at Virginia Tech.

\bibliographystyle{plain}
\bibliography{references}
\end{document}